# Inverse problem for one class of nonselfadjoint operator's bunches with nonperiodic coefficients

Institute of Applied Mathematics, Baku State University, Azerbaijan
370148, Baku, Z.Khalilov, 23,
rakibaz@yahoo.com
R.F.Efendiev

**Abstract:**

In this paper the complete spectral analysis of the operators is carried out and also with help of generalized normalizing numbers the inverse problem is solved.

The present work is devoted to studying of a differential operator's bunch $L(k)$ generated by differential expression

$$l(y) \equiv (-1)^m y^{(2m)}(x) + \sum_{\gamma=0}^{2m-2} p_\gamma(x,k) y^{(\gamma)}(x) - k^{2m} y \tag{1}$$

in the space $L_2[0,\infty)$, with boundary conditions

$$y^{(j)}(0) = 0, \; j = \overline{0, m-1} \tag{2}$$

under the assumption, that

$$p_\gamma(x,k) = \sum_{s=0}^{2m-\gamma-1} \sum_{n=1}^{\infty} p_{\gamma s n} k^s e^{-nx}, \; 0 \le x < \infty \tag{3}$$

and a series $\sum_{\gamma=0}^{2m-2} \sum_{s=0}^{2m-\gamma-1} \sum_{n=1}^{\infty} n^{\gamma+s} |p_{\gamma s n}| \tag{4}$

converges. In the work the data of scattering of the boundary problem (1) - (2) are entered and the theorem of uniqueness of definition $L(k)$ on them is proved. We shall note, that at $p_\gamma(x,k) = \sum_{n=1}^{\infty} p_{\gamma 0 n} e^{-nx}$ these problem has been considered by M.G.Gasymov [1].

Special solution $l(y) = 0$.

Let's introduce the following designations:

$\omega_j = \exp(ij\pi/m), \; j = 0,1,...2m-1$

$k_{nj\tau} = -in/\omega_\tau (1-\omega_j), \; n = 1,2,....; \; \tau = 0,1,...2m-1;$

$j = 1,2,...2m-1; \; k_{nj} = k_{nj0}$

$C_m^n = \dfrac{m!}{(m-n)! n!}$

$\dfrac{1}{in + k\omega_\tau(1-\omega_j)} \left[ (i\alpha + k\omega_\tau)^{2m} - k^{2m} - (i\alpha + k_{nj})^{2m} - k_{nj}^{2m} \right] =$

$= \sum_{\gamma=0}^{2m-2} C_{j\tau\gamma}(n,\alpha) k^\gamma \cdot \tau = 0,1,...2m-1; \; j = 1,2,...2m-1; n = 1,2...$

$$\frac{1}{in + k\omega_\tau (1-\omega_j)} \left[ k^s (ik\omega_\tau - t)^\gamma - k^s_{nj\tau} (ik_{nj} - t)^\gamma \right] = \sum_{v=0}^{\gamma-1} C^{(n,r,v)}_{j\tau vs} k^{v+s}$$

**Theorem 1.** Suppose that condition (2-4) satisfied. Then the differential equation $l(y) = 0$ has the special solutions $f_\tau(x,k)$, $\tau = \overline{0,2m-1}$, representable as

$$f_\tau(x,k) = e^{ik\omega_\tau x} + \sum_{\alpha=1}^{\infty} V_\alpha^{(\tau)} e^{(-\alpha + ik\omega_\tau)x} + \sum_{j=0}^{2m-1} \sum_{\alpha=1}^{\infty} \sum_{n=1}^{\alpha} \frac{V_{n\alpha}^{(j,\tau)}}{in + k\omega_\tau (1-\omega_j)} e^{(-\alpha + ik\omega_\tau)x} \quad (5)$$

where numbers $V_{n\alpha}^{(j,\tau)}$, $V_\alpha^{(\tau)}$ are determined from the following recurrent formulas

$$\left[ (i\alpha + k_{nj})^{2m} - k^s_{nj} \right] V_{n\alpha}^{(j,\tau)} + + \sum_{\gamma=0}^{2m-2} \sum_{s=0}^{2m-\gamma-1} \sum_{r=n}^{\alpha-1} k^s_{nj\tau} (ik_{nj} - r)^\gamma P_{\gamma s\alpha - r} V_{nr}^{(j,\tau)} = 0 \quad (6)$$

at $\alpha = 2,3,\ldots; n = 1,2,\ldots\alpha-1; j = 1,\ldots 2m-1; \tau = 0,\ldots 2m-1;$

$$C^\gamma_{2m}(i\alpha)^{2m-\gamma} \omega^\gamma_\tau V_\alpha^{(\tau)} + \sum_{s=0}^{\gamma} P_{s,\gamma-s,\alpha}(i\omega_\tau)^s + \sum_{j=1}^{2m-1} \sum_{n=1}^{\alpha} C_{j\tau\gamma}(n,\alpha) V_{n\alpha}^{(j,\tau)} +$$

$$+ \sum_{s=0}^{\gamma} \sum_{v=s}^{2m-2} \sum_{r=1}^{\alpha-1} (-1)^{v-s} C^s_v (i,\omega_\tau)^s (ir)^{v-s} P_{v,\gamma-s,\alpha-2} V_r^{(\tau)} +$$

$$+ \sum_{j=1}^{2m-1} \sum_{v=0}^{\gamma} \sum_{s=v+1}^{2m-2} \sum_{r=n}^{\alpha-1} C_{j,v,\tau,\gamma-v}(n,t,v) P_{s,\gamma-v,\alpha-2} V_{nr}^{(j,\tau)} = 0 \quad (7)$$

at $\gamma = \overline{0,2m-2}; \tau = \overline{0,2m-2}; \alpha = 1,2,\ldots, n < \alpha$.

$$2mi\alpha\omega_\tau^{2m-1} V_\alpha^{(\tau)} + \sum_{s=0}^{2m-1} P_{s,2m-1-s,\alpha}(i\omega_\tau)^s + \sum_{s=0}^{2m-2} \sum_{r=1}^{\alpha-1}(i\omega_\tau) P_{s,2m-1-s,r} V_{\alpha-r}^{(\tau)} = 0 \quad (8)$$

and series

$$\sum_{j=1}^{2m-1} \sum_{n=1}^{\infty} \frac{1}{n} \sum_{\alpha=n}^{\infty} \alpha^{2m-1}(\alpha-n) \left| V_{n\alpha}^{(j,\tau)} \right|,$$

$$\sum_{j=1}^{2m-1} \sum_{\alpha=1}^{\infty} \alpha^{2m-1} \left| V_{\alpha\alpha}^{(j,\tau)} \right|,$$

$$\sum_{\alpha=1}^{\infty} \alpha^{2m} \left| V_\alpha^{(\tau)} \right|,$$

converge.

**Theorem 2.** At any values of $n$ and $j$ $\nu = \overline{0,2m-1}$ $j = 1,2,\ldots 2m-1$, $n = 1,2\ldots;$ has places the equality

$$f^\nu_{nj}(x) = V_{nn}^{(j,\nu)} f_{j+\nu}(x, k_{nj\nu}) \quad (9)$$

**where**

$$f^\nu_{nj}(x) = \operatorname{Res} f_\nu(x,k) \big|_{k=k_{nj\nu}}$$

$$f_{2m+\nu}(x,k) = f_\nu(x,k)$$

If to break $K$ - a complex plane, into sectors $S_\nu = \left\{ \frac{\nu\pi}{2} < \arg k < \frac{(\nu+1)\pi}{2}, \nu = \overline{0,2m-1} \right\}$

and to number $f_\nu(x,k)$ so that in sector $S_\nu$ it was fulfilled

$$\operatorname{Im}(k\omega_0) < \operatorname{Im}(k\omega_1) < \ldots < \operatorname{Im}(k\omega_{m-1}) < 0 < \operatorname{Im}(k\omega_m) < \ldots < \operatorname{Im}(k\omega_{2m-1}) \quad (10)$$

then it is obvious, that $f_\nu(x,k) \in L_2(0,\infty)$ at $\nu = \overline{0,m-1}$, $f_\nu(x,k) \overline{\in} L_2(0,\infty)$

$\tau = \overline{m, 2m-1}$

**Theorem 3.** In order to $k$, $\text{Im } k \neq 0$ be eigenvalue of the bunch $L(k)$ it is necessary and sufficient the equation is fulfilled

$$W_v(k) = \begin{vmatrix} f_0(0,k) & f_1(0,k) \ldots \ldots f_{m-1}(0,k) \\ \vdots & \\ f_0^{(m-1)} & f_1^{(m-1)}(0,k) \ldots f_{m-1}^{(m-1)}(0,k) \end{vmatrix} = 0, k \in S_v \quad (11)$$

**Theorem 4.** The bunch $L(k)$ has final number of eigenvalue on each of sectors $S_v, v = \overline{m, 2m-1}$, included boundary.

II. Constructions of a kernel resolvent of the bunch $L(k)$.

Let's calculate the kernel of the resolvent operator's bunch $L(k)$. .

Let $\psi(x) \in L_2[0, \infty)$. Applying the variation method, we search for the solution of the equation $l(y) = \varphi(x)$ from $L_2[0, \infty)$ as

$$y(x,k) = \sum_{i=0}^{2m-1} C_i(x,k) f_i(x,k).$$

Then we obtain the following system of the equations.

$$\sum_{i=0}^{2m-1} C_i'(x,k) f_i^{(j)}(x,k) = 0, \text{ at } j = \overline{0, 2m-2}$$

$$\sum_{i=0}^{2m-1} C_i'(x,k) f_i^{(2m-1)}(x,k) = \psi(x)$$

Solving them we get

$$C_s'(x,k) = -\frac{i\omega_s}{4k^{2m-1}} \varphi_s(x,k)\psi(x), \ s = \overline{0, 2m-1}$$

**Theorem 5.** Functions $\varphi_s(x,k)$ $s = \overline{0, 2m-1}$ have the following properties

$$\varphi_s(x,k) = e^{-ik\omega_s x} + \sum_{\alpha=1}^{\infty} R_\alpha^{(s)} e^{(-\alpha - ik\omega_s)x} + \sum_{j=0}^{2m-1}\sum_{n=1}^{\infty}\sum_{\alpha=n}^{\infty} \frac{R_{n\alpha}^{(j,s)}}{in - k\omega_s(1-\omega_j)} e^{(-\alpha - ik\omega_s)x} \quad (12)$$

Functions $\varphi_s(x,k)$ $s = \overline{0, 2m-1}$ are solutions of the equation.

$$(-1)^m Z^{(2m)}(x) + \sum_{\gamma=0}^{2m-2} (-1)^\gamma \left[P_\gamma(x,k)Z\right]^{(\gamma)} = k^{2m} Z$$

Remark 1. Taking into account (10) we shall obtain $\varphi_s(x,k) \in L_2[0, \infty)$ at $s = \overline{m, 2m-1}$, $\varphi_s(x,k) \overline{\in} L_2[0, \infty)$ at $s = \overline{0, m-1}$.

Let's continue calculation of the kernel of the resolvent, considering the property of functions $f_s(x,k), \varphi_s(x,k)$, and also, that $y(x,k) \in L_2[0, \infty)$

$$C_s(x,k) = -\frac{i\omega_s}{2mk^{2m-1}} \int_0^x \varphi_s(t,k)\psi(t)dt + C_s(k). \ s = \overline{0, 2m-1}.$$

Then

$$y(x,k) = \sum_{s=0}^{2m-1} \left[ -\frac{i\omega_s}{2mk^{2m-1}} \int_0^x f_s(x,k)\varphi_s(t,k)\psi(t)dt + C_s(k)f_s(x,k). \right]$$

Now let's construct the kernel of the resolvent on sector $S_v$. Taking into account, that $\varphi_s(x,k) \in L_2[0, \infty)$ at $s = \overline{m, 2m-1}$, we get

$$C_s(k) = \frac{i\omega_s}{2mk^{2m-1}} \int_0^\infty \varphi_s(t,k)\psi(t)dt; \quad s = \overline{m, 2m-1} \qquad (13)$$

Then

$$y(x,k) = \sum_{s=0}^{m-1} \left[ f_s(x,k) \left[ C_s(k) - \frac{i\omega_s}{2mk^{2m-1}} \int_0^x \varphi_s(t,k)\psi(t)dt \right] \right] +$$

$$+ \sum_{s=0}^{2m-1} \frac{i\omega_s}{2mk^{2m-1}} f_s(x,k) \int_x^\infty \varphi_s(t,k)\psi(t)dt$$

As $y(x,k)$ should satisfy to the condition $y^{(j)}(0) = 0, \ j = \overline{0, m-1},$ for definition $C_s(k), \ s = \overline{0, m-1}$ we should get the following system of the equations

$$\sum_{s=0}^{m-1} f_s^{(j)}(0,k) C_s(k) = -\sum_{s=m}^{2m-1} f_s^{(j)}(0,k); \quad j = \overline{0, m-1}.$$

Let's designate through

$$A_{j,s}(k) = \begin{vmatrix} f_0(0,k), \ldots f_{j-1}(0,k) f_s(0,k) \cdot f_{j+k}(0,k) \ldots \ldots f_{m-1}(0,k) \\ \ldots \ldots \\ f_0^{(m-1)}(0,k) \ldots f_{j-1}^{(m-1)}(0,k) f_s^{(m-1)}(0,k) f_{j+1}^{(m-1)}(0,k) \ldots f_{m-1}^{(m-1)}(0,k) \end{vmatrix}$$

$j = \overline{0, m-1}, \ s = \overline{m, 2m-1}.$

Then for $C_j(k), \ j = \overline{0, m-1}$ we have

$$C_j(k) = \sum_{s=m}^{2m-1} C_s(k) \frac{A_{js}(k)}{W_v(k)};$$

Considering (13) we get

$$C_j(k) = \sum_{s=m}^{2m-1} \frac{i\omega_s}{2mk^{2m-1}} \int_0^\infty \varphi_s(t,k)\psi(t)dt \cdot \frac{A_{js}(k)}{W_v(k)}; \quad j = \overline{0, m-1}$$

Then

$$y(x,k) = \sum_{j=0}^{m-1} f_j(x,k) \left[ \sum_{s=m}^{2m-1} \frac{i\omega_s}{2mk^{2m-1}} \int_0^\infty \varphi_s(t,k)\psi(t)dt \cdot \frac{A_{js}(k)}{W_v(k)} - \right.$$

$$\left. - \frac{i\omega_j}{2mk^{2m-1}} \int_0^\infty \varphi_j(t,k)\psi(t)dt \right] + \sum_{j=m}^{2m-1} \frac{i\omega_j}{2mk^{2m-1}} f_j(x,k) \int_x^\infty \varphi_j(t,k)\psi(t)dt$$

Then for the kernel $R_v(x,t,k), \ k \in S_v$ we have

$$R_v(x,t,k) = \frac{1}{2mk^{2m-1}} \begin{cases} \sum_{j=0}^{m-1} \left[ f_j(x,k) \sum_{s=m}^{2m-1} i\omega_s \varphi_s(t,k) \frac{A_{is}(k)}{W_v(k)} - \right. \\ \left. - i\omega_j \varphi_j(t,k) \right], \ npu \ t < x \\ \sum_{j=0}^{m-1} \sum_{s=m}^{2m-1} i\omega_s f_j(x,k) \varphi_s(t,k) \frac{A_{is}(k)}{W_v(k)} + \\ + \sum_{j=m}^{2m-1} i\omega_j f_j(x,k) \varphi_j(t,k), \ npu \ t > x \end{cases}$$

at $k \in S_v$

Now we shall find a kind of decomposition on eigenfunctions of a bunch $L(k)$. Let $f(x)$ - finite function, the support of which does not involve the point $x = 0$. Using equality

$$(-1)^m \frac{\partial^{2m}}{\partial x^{2m}} R(x,t,k) + \sum_{\gamma=0}^{2m-2} P_\gamma(x,k) \frac{\partial^\gamma}{\partial x^\gamma}[R(x,t,k)] - k^{2m} R(x,t,k) = \delta(x,t)$$

and integrating in parts, we receive, that

$$\int_0^\infty R(x,t,k) f(t) dt = \frac{f(x)}{k^{2m}} + 0\left(\frac{1}{k^{2m-1}}\right).$$

Remark. As it has above been proved on each sector $S_\nu$, including their borders, there is a final number of eigenvalues. We shall designate them through $k_1^\nu,....,k_{j\nu}^\nu$ on sector $S_\nu, \nu = \overline{0,2m-1}$. Let $\lambda_{1,\nu},....\lambda_{\tau,\nu}$, $\tau < j_\nu$ be those, points where some from numbers coincide if such numbers exist and $\lambda_{1,\nu} < .... < \lambda_{\tau,\nu}$. Then through $C_\gamma$ we shall designate a contour
$(-\infty, \lambda_{1,\nu} - \delta), (\lambda_{1,\nu} + \delta, \lambda_{2,\nu} - \delta),...., (\lambda_{\tau,\nu} + \delta, \lambda_{\tau,\nu} - \delta)$ semicircles $\delta$ with the centers in points $\lambda_{1,\nu},.....\lambda_{\tau,\nu}$. Here $\delta < \varepsilon/2$.

Let $\Gamma_N$ - be a contour of a circle of radius $N$ with the center in zero. Then at $N \to \infty$ equalities are true

$$\frac{1}{2\pi i} \oint_{\Gamma_N} k^i y(x,k) dk = 0(1), \quad i = \overline{0, 2m-2}$$

$$\frac{1}{2\pi i} \oint_{\Gamma_N} k^{2m-1} y(x,k) dk = f(x) + 0(1)$$

and under theorem Cauchy.

$$\frac{1}{2\pi i} \oint_{\Gamma_N} k^{2m-1} y(x,k) dk = \frac{1}{2\pi i} \int_{\Gamma_N} k^{2m-1} R(x,t,k) f(t) dt dk =$$

$$= \frac{1}{2\pi i} \sum_{\nu=0}^{2m-1} \int_{C_\nu} \int_0^\infty k^{2m-1} [R_\nu(x,t,k) - R_{\nu-1}(x,t,k)] f(t) dt dk -$$

$$- \sum_{\nu=0}^{2m-1} \int_0^\infty \operatorname{Res} R_\nu(x,t,k) k^{2m-1} \bigg|_{k=\lambda_{\tau,\nu}} \cdot f(t) dt = f(x)$$

$$= \frac{1}{2\pi i} \sum_{\nu=0}^{2m-1} \oint_{C_\nu} \int_0^\infty k^i [R_\nu(x,t,k) - R_{\nu-1}(x,t,k)] f(t) dt dk -$$

$$- \sum_{\nu=0}^{2m-1} \int_0^\infty \operatorname{Res} R_\nu(x,t,k) k^i \bigg|_{k=\lambda_{\tau,\nu}} \cdot f(t) dt = 0, \quad i = \overline{0, 2m-2}$$

**Inverse problem:**

Set of value $S_j(k) = \frac{A_{ij}(k)}{W_\nu(k)}, j = \overline{m, 2m-1}$, $S_j^{-1}(k)$ we shall name the scattering data of the boundary problem (1) - (2). The inverse problem of definition of an operator's bunch $L(k)$ in them is put.

Here is the algorithm of construction of the coefficients $P_\gamma(x,k), \gamma = \overline{0, 2m-2}$.

1. On the value $S_j(k)$ and $S_j^{-1}(k)$ all numbers are determined $V_{nn}^{(j,\nu)}$;

$j = 1, 2m-1, \ v = 0, 2m-1, n = 1, 2, \ldots$

2. Under the formula

$$V_{m,\alpha+m}^{(j,v)} = V_{m,m}^{(j,v)} (1-\omega_j) \left( V_{\alpha}^{(j+v)} + \sum_{s=1}^{2m-1} \sum_{r=1}^{\alpha} \frac{V_{r\alpha}^{(s,j,v)}}{r(1-\omega_j) - n\omega_j(1-\omega_s)} \right)$$

at $j = 1,2,\ldots 2m-1, \ \alpha = 1,2,\ldots, n = 1,2,\ldots \alpha - 1$; all numbers are determined $V_{n\alpha}^{(j,v)}, V_{\alpha}^{(v)}$

3. Using equalities (6-8) is all numbers $P_{\gamma sn}$,

$\gamma = 0, 2m-2, \ s = 0, 2m - \gamma - 1, \ n = 1, 2 \ldots$

In more detail we shall consider definition of numbers $V_{nn}^{(j,v)}$.

$$V_{nn}^{(j,s)} = \lim_{k \to k_{n,i-s,s}} (k - k_{n,i-s,s}) \frac{A_{is}(k)}{W_v(k)}, \qquad s = \overline{m, 2m-1}$$

$$V_{nn}^{(j,s)} = \lim_{k \to k_{n,i-s,s}} (k - k_{n,i-s,s}) \frac{W_v(k)}{A_{is}(k)}, \qquad s = \overline{0, m-1}$$

As numbers $V_{n\alpha}^{(j,v)} \ \alpha = 1,2,\ldots n = 1,2\ldots, n < \alpha$ are determined unequivocal and is unique it is possible to formulate the theorem of uniqueness.

Theorem 6. Let set of values $S_j(k), \ S_j^{-1}(k), \ j = \overline{m, 2m-1},$, which is given, is by the scattering data of some operator's bunch of type $L(k)$ at performance of conditions (2) - (4). Then the operator's bunch $L(k)$ answering to them is uniqueness.

Reference:
1. Casymov M.G. Uniqueness of the solutions of an inverse problem of the theory of scattering for one class of ordinary differential operators even orders. Report AS.USSR. 1982,V. 266,N5, p 1033-1036.